\documentclass{compositio}

\theoremstyle{plain}
\newtheorem{theorem}{Theorem}[section]
\newtheorem{lemma}{Lemma}[section]
\newtheorem{corollary}{Corollary}[section]
\newtheorem{proposition}{Proposition}[section]

\theoremstyle{definition} 
\newtheorem{definition}[theorem]{Definition} 
\newtheorem{conjecture}{Conjecture}

\theoremstyle{remark} 
\newtheorem{remark}{Remark}
\newtheorem{example}{Example}

\begin{document}

\title{Non-rational nodal quartic threefolds}
\author{Ivan Cheltsov}
\email{cheltsov@yahoo.com}
\address{Steklov Institute of Mathematics\hfill\break\indent 8 Gubkin street, Moscow 117966\hfill\break\indent Russia\hfill\break\indent }
\thanks{All varieties are assumed to be projective, normal and defined over ${\mathbb C}$.}%

\begin{abstract}
The ${\mathbb Q}$-factoriality of a nodal quartic 3-fold implies
its non-rationality. We prove that a nodal quartic 3-fold with at
most $8$ nodes is ${\mathbb Q}$-factorial, and we show that a
nodal quartic 3-fold with $9$ nodes is not ${\mathbb Q}$-factorial
if and only if it contains a plane. However, there are
non-rational non-${\mathbb Q}$-factorial nodal quartic 3-folds in
${\mathbb P}^4$. In particular, we prove the non-rationality of a
general non-${\mathbb Q}$-factorial nodal quartic 3-fold that
contains either a plane or a smooth del Pezzo surface of degree
$4$.
\end{abstract}

\maketitle


\section{Introduction}
\label{section:introduction}

Consider a nodal quartic 3-fold $X\subset{\mathbb P}^{4}$, i.e. a
hypersurface of degree $4$ whose singular points are simple double
points. The following result was proved in \cite{IsMa71},
\cite{Pu88b} and \cite{Me03}.

\begin{theorem}
\label{theorem:of-Mella} Let $X$ be ${\mathbb Q}$-factorial. Then
$X$ is not birational to a ${\mathbb Q}$-factorial terminal Fano
3-fold with Picard group ${\mathbb Z}$ that is not biregular to
$X$, and $X$ is not birational to a fibration of varieties of
Kodaira dimension $-\infty$.
\end{theorem}

In this paper we will prove the following result.

\begin{theorem}
\label{theorem:main}%
Suppose $|\textsf{Sing}(X)|\le 8$. Then $X$ is ${\mathbb
Q}$-factorial.
\end{theorem}

\begin{corollary}
\label{corollary:main} Nodal quartic 3-folds with at most $8$
nodes are non-rational.
\end{corollary}

The conditions of Theorem~\ref{theorem:main} can not be weakened.

\begin{example}
\label{example:quartic-plane} Let $X$ be a sufficiently general
quartic 3-fold containing a two-dimensional linear subspace
$\Pi\subset{\mathbb P}^4$. Then $X$ is nodal and non-${\mathbb
Q}$-factorial, the quartic $X$ has $9$ nodes, which are the
intersection of two cubic curves in the plane $\Pi$.
\end{example}

However, we will prove the following result.

\begin{theorem}
\label{theorem:second}%
Suppose $|\textsf{Sing}(X)|=9$. Then $X$ is non-${\mathbb
Q}$-factorial if and only if it contains a plane.
\end{theorem}

\begin{remark}
\label{remark:nine-points}%
A general nodal quartic 3-fold with $9$ nodes is ${\mathbb
Q}$-factorial (see \cite{GrSt} and \cite{EiGr01}).
\end{remark}

A posteriori the non-${\mathbb Q}$-factoriality of the quartic $X$
does not necessarily imply its rationality, viz. we will prove
the following result (cf. Remark 3 in \cite{Me03}). 

\begin{theorem}
\label{theorem:quartic-plane} Let $X$ be a very general quartic
3-fold containing a plane. Then $X$ is non-rational.
\end{theorem}

Rational nodal quartic 3-folds do exist.

\begin{example}
\label{example:determinantal-quartic} Let $X$ be a general
determinantal quartic 3-fold. Then $X$ is nodal, non-${\mathbb
Q}$-factorial and rational, and $|{\textsf{Sing}}(X)|=20$ (see
\cite{Pet98} and \cite{Me03}).
\end{example}

\begin{remark}
\label{remark:number-of-points}%
The quartic $X$ can not have more than $45$ nodes by \cite{Va83}
and \cite{Fr86}, and $X$ can have any number of nodes up to $45$
(see \cite{CiGe03}). There is a unique (see \cite{JSV90}) nodal
quartic 3-fold ${\mathcal B}_{4}$ with $45$ nodes which can be
given by the equation
$$
w^4-w(x^3+y^3+z^3+t^3)+3xyzt=0\subset {\mathbb P}^4\cong{\textsf{Proj}}({\mathbb C}[x,y,z,t,w])%
$$
and is known as the Burkhardt quartic (see \cite{Bu90},
\cite{Bu92}, \cite{To36}, \cite{Ba46}, \cite{Fi89}, \cite{Pet98}).
The quartic ${\mathcal B}_{4}$ is determinantal. Moreover, the
quartic ${\mathcal B}_{4}$ is the unique invariant of degree $4$
of the simple group ${\textsf{PSp}}(4,{\mathbb Z}_{3})$ of order
$25920$ (see \cite{vdGr87}, \cite{Hu96}, \cite{HoWe01} and
\cite{HulSa02}). The nodes of ${\mathcal B}_{4}$ correspond to the
$45$ tritangents of a smooth cubic surface, and the Weyl group of
$E_{6}$ is a nontrivial extension of ${\textsf{PSp}}(4,{\mathbb
Z}_{3})$ by ${\mathbb Z}_{2}$.
\end{remark}

For a given variety, it is the one of the most substantial
questions to decide whether it is rational or not. This question
was considered in depth for smooth 3-folds (see \cite{Ro55},
\cite{IsMa71}, \cite{ClGr72}, \cite{Be77}, \cite{Tu80},
\cite{Ti80a}, \cite{Ti80b}, \cite{Is80b}, \cite{Sar80},
\cite{Ti82}, \cite{Sar82}, \cite{Sho84}, \cite{Ti86}, \cite{Al87},
\cite{Pu89}, \cite{Co95}, \cite{Is97}, \cite{Pu98b},
\cite{IsPr99}, \cite{Co00}). On the other hand, even relatively
mild singularities can force a 3-fold to be rational. For example,
with a few exceptions all canonical Gorenstein Fano 3-folds having
a non-cDV point are rational due to \cite{Pr04}, but in the
non-Gorenstein case the situation is different (see \cite{CPR},
\cite{Ch97}, \cite{Ch04}). Therefore, the rationality of nodal
3-folds can be considered as a rather natural topic (see
\cite{Pu88b}, \cite{Gr88a}, \cite{Gr88b}, \cite{Me03},
\cite{ChPa04}).

\begin{remark}
\label{remark:non-rational-threefolds} Every nodal hypersurface in
${\mathbb P}^4$ of degree at least $5$ is non-rational. All
quadric 3-folds are rational. A nodal cubic 3-fold in ${\mathbb
P}^4$ is non-rational if and only if it is smooth by
\cite{ClGr72}.
\end{remark}

There are non-${\mathbb Q}$-factorial nodal quartic 3-folds that
contain no planes (see \cite{Me03} and \cite{El00}).

\begin{example}
\label{example:quartic-quadric} Consider a sufficiently general
quartic 3-fold $X\subset{\mathbb P}^4$ passing through a smooth
quadric surface $Q\subset{\mathbb P}^4$. The quartic $X$ can be
given by the equation
$$
a_{2}(x,y,z,t,w)h_{2}(x,y,z,t,w)=b_{3}(x,y,z,t,w)g_{1}(x,y,z,t,w)\subset{\textsf{Proj}}({\mathbb C}[x,y,z,t,w]),%
$$
where $a_{2}$, $h_{2}$, $b_{3}$ and $g_{1}$ are homogeneous
polynomials of degree $2$, $2$, $3$ and $1$ respectively, and the
quadric surface $Q\subset{\mathbb P}^4$ is given by the equations
$h_{2}=g_{1}=0$. The quartic $X$ is non-${\mathbb Q}$-factorial
and nodal, and it has $12$ nodes given by
$h_{2}=g_{1}=a_{2}=b_{3}=0$. Introducing a new variable
$\alpha=a_{2}/g_{1}$ one can unproject $X$ (see \cite{Re00}) into
a complete intersection $V\subset{\mathbb P}^5$ given by the
equations
$$
\alpha g_{1}(x,y,z,t,w)-a_{2}(x,y,z,t,w)=\alpha h_{2}(x,y,z,t,w)-b_{3}(x,y,z,t,w)=0\subset {\mathbb P}^{5}%
$$
such that the unprojection $\rho:X\dasharrow V$ is a composition
$\rho=\phi\circ\psi^{-1}$, where $\phi:Y\to V$ is an extremal
contraction and $\psi:Y\to X$ is a flopping contraction (see
\cite{Ko89}, \cite{Co95}). The variety $V$ is smooth outside of a
point $P=(0:0:0:0:0:1)$ which is a node on $V$. The morphism
$\phi$ contracts the surface ${\mathbb P}^1\times {\mathbb P}^1$
to $P$, and $\psi$ contracts the images of $12$ lines on $V$
passing through $P$ into the nodes of $X$. It is unknown whether
$X$ is rational or not (see \cite{Is80b}, \cite{Pu89},
\cite{IsPu96} and \cite{Co00}).
\end{example}

There are non-rational non-${\mathbb Q}$-factorial nodal quartic
3-folds in ${\mathbb P}^4$ that do not contain neither planes nor
quadric surfaces. In particular, we will prove the following
result.

\begin{theorem}
\label{theorem:quartic-del-Pezzo} Let $X\subset{\mathbb P}^4$ be a
sufficiently general quartic 3-fold containing a smooth del Pezzo
surface $S\subset{\mathbb P}^4$ of degree $4$. Then $X$ is nodal,
non-${\mathbb Q}$-factorial and non-rational,
$|{\textsf{Sing}}(X)|=16$.
\end{theorem}

The quartic 3-folds in Theorems~\ref{theorem:quartic-plane} and
\ref{theorem:quartic-del-Pezzo} are the only known examples of
nodal, non-rational and non-${\mathbb Q}$-factorial quartic
3-folds. The degeneration technique (see \cite{Be77}, \cite{Tu80},
\cite{Cl83}, \cite{Ko95}, \cite{Ko96}) together with either
Theorem~\ref{theorem:quartic-plane} or
Theorem~\ref{theorem:quartic-del-Pezzo} give another proof that a
very general smooth quartic 3-fold is non-rational (see
\cite{IsMa71}, \cite{CoMu77a} and \cite{CoMu77b}).

\begin{remark}
\label{remark:Luroth} There are few known examples of unirational
smooth quartic 3-folds (see \cite{Seg60}, \cite{IsMa71},
\cite{Is80b} and \cite{Mar00}). Moreover, it is still unknown
whether a generic smooth quartic 3-fold is unirational or not.
However, the quartics in Theorems~\ref{theorem:quartic-plane} and
\ref{theorem:quartic-del-Pezzo} are birational to fibrations of
del Pezzo surfaces of degree $3$ or $4$. Thus, the quartics in
Theorems~\ref{theorem:quartic-plane} and
\ref{theorem:quartic-del-Pezzo} are unirational (see \cite{Ma66},
\cite{Ma67}, \cite{Ma72}).
\end{remark}

Both Theorems~\ref{theorem:main} and \ref{theorem:second} can be
considered as a part of the following conjecture.

\begin{conjecture}
\label{conjecture:factorial-quartics} Let $V\subset{\mathbb P}^4$
be a nodal hypersurface. Then $V$ is ${\mathbb Q}$-factorial if
one of the following three conditions holds:
\begin{itemize}
\item $|{\textsf{Sing}}(V)|<(\deg(V)-1)^2$;%
\item $|{\textsf{Sing}}(V)|<2(\deg(V)-1)(\deg(V)-2)$ and $V$ contains no planes;%
\item $|{\textsf{Sing}}(V)|\le 2(\deg(V)-1)(\deg(V)-2)$ and $V$ contains neither planes nor quadrics.%
\end{itemize}
\end{conjecture}

It is easy to see that
Conjecture~\ref{conjecture:factorial-quartics} holds for quadrics
and cubics (see \cite{FiWe89}). Moreover, an analogue of
Conjecture~\ref{conjecture:factorial-quartics} for smooth surfaces
on a nodal hypersurface in ${\mathbb P}^4$ is proven in
\cite{CiGe03}.

\begin{acknowledgements}
The author is grateful to A.~Corti, M.~Grinenko, V.~Iskovskikh,
S.~Kudryavtsev, M.~Mel\-la, J.~Park, Yu.~Pro\-kho\-rov,
A.~Pan\-kra\-tiev, A.~Pukh\-li\-kov, V.~Sho\-ku\-rov and
L.~Wotzlaw for fruitful conversations.
\end{acknowledgements}

\section{The proof of Theorems~\ref{theorem:main} and \ref{theorem:second}}
\label{section:proof-of-main-theorems}

The ${\mathbb Q}$-factoriality of a nodal 3-fold depends on the
global position of its nodes. This subject was thoroughly studied
in \cite{Cl83}, \cite{Sch85}, \cite{Fi87}, \cite{We87},
\cite{FiWe89}, \cite{Di90}, \cite{Bo90}, \cite{Cy99}, \cite{En99},
\cite{Cy01}. In particular, let $X\subset{\mathbb P}^{4}$ be a
nodal quartic 3-fold.

\begin{remark}
\label{remark:factroiality} The following conditions are
equivalent (see \cite{We87}, \cite{Di90}, \cite{Cy01},
\cite{CoMe02}):
\begin{itemize}
\item the quartic $X$ is ${\mathbb Q}$-factorial;%
\item every Zariski local ring of the quartic $X$ is UFD, i.e. $X$ is factorial;%
\item the group $H_{4}(X,{\mathbb Z})$ is generated by the class of a hyperplane section;%
\item $\dim(H_{4}(X,{\mathbb Z}))=\dim(H^{2}(X,{\mathbb Z}))=1$;%
\item the nodes of $X$ impose independent linear conditions on
cubic hypersurfaces in ${\mathbb P}^4$.%
\end{itemize}
\end{remark}

Suppose that $X$ does not contain planes and
$|{\textsf{Sing}}(X)|\le 9$. We will show that the nodes of the
quartic $X$ impose independent linear conditions on cubic
hypersurfaces in ${\mathbb P}^4$.

\begin{definition}
\label{definition-of-general-position} The points of a set
$\Gamma\subset{\mathbb P}^4$ are in almost general position if the
following holds:
\begin{itemize}
\item at most $3$ points of $\Gamma$ can lie on a line;%
\item at most $6$ points of $\Gamma$ can lie on a conic;%
\item at most $8$ points of $\Gamma$ can lie on a plane;%
\end{itemize}
\end{definition}

\begin{proposition}
\label{proposition:nodes-in-general-position} The nodes of the
quartic $X$ are in almost general position.
\end{proposition}

\begin{proof}
Let $L\subset{\mathbb P}^4$ be a line and $\Pi\subset{\mathbb
P}^4$ be a sufficiently general two-dimensional linear subspace
passing through $L$. Then $\Pi\not\subset X$ and $\Pi\cap X=L\cup
S$, where $S$ is a plane cubic curve. Moreover,
$$
{\textsf{Sing}}(X)\cap L\subset L\cap S,
$$
but $|L\cap S|\le 3$. Thus, at most $3$ nodes of the quartic $X$
can lie on a line in ${\mathbb P}^4$.

Let $C\subset{\mathbb P}^4$ be a smooth conic and
$Y\subset{\mathbb P}^4$ be a sufficiently general two-dimensional
quadric cone over $C$. Then $Y\not\subset X$ and $Y\cap X=C\cup
R$, where $R$ is a curve of degree $6$. As above we have the
inclusion ${\textsf{Sing}}(X)\cap C\subset C\cap R$. However, the
curves $C$ and $R$ lie in the smooth locus of $Y$ and the
intersection $C\cdot R$ on $Y$ equals to $6$. Thus, the inequality
$|C\cap R|\le 6$ holds. Hence, at most $6$ nodes of the quartic
$X$ can lie on a smooth conic in ${\mathbb P}^4$.

Let $\Sigma\subset{\mathbb P}^4$ be a plane and $T=\Sigma\cap X$.
Then $T$ is a possibly reducible and non-reduced plane quartic and
${\textsf{Sing}}(X)\cap\Sigma\subset{\textsf{Sing}}(T)$. In
particular, $|{\textsf{Sing}}(X)\cap\Sigma|\le 6$ in the case of
non-reduced curve $T$, because we already proved that at most $3$
nodes of $X$ can lie on a line and at most $6$ nodes of the
quartic $X$ can lie on a conic. However, $|{\textsf{Sing}}(T)|\le
6$ whenever $T$ is reduced. Therefore, at most $6$ nodes of $X$
can lie on a plane in ${\mathbb P}^4$.
\end{proof}

\begin{proposition}
\label{proposition:nodes-plane} Let $\Pi\subset{\mathbb P}^4$ be a
two-dimensional linear subspace such that ${\textsf{Sing}}(X)$ is
contained in $\Pi$. Then the nodes of $X$ impose independent
linear conditions on cubic curves in $\Pi\cong{\mathbb P}^2$ and
on cubic hypersurfaces in ${\mathbb P}^4$.
\end{proposition}

\begin{proof}
We must show that for any subset
$\Sigma\subsetneq{\textsf{Sing}}(X)$ and a point
$p\in{\textsf{Sing}}(X)\setminus\Sigma$ there is a cubic curve in
$\Pi$ and a cubic hypersurface in ${\mathbb P}^4$ passing through
the points in $\Sigma$ and not passing through the point $p$. Let
$\pi:V\to\Pi$ be a blow up of points in $\Sigma$. Then $V$ is a
weak del Pezzo surface of degree $9-|\Sigma|\ge 2$ due to
Proposition~\ref{proposition:nodes-in-general-position} and
$|-K_{V}|$ is free (see \cite{De80}, \cite{HiWa81}, \cite{Bes83},
\cite{Mae94}). There is a curve $C\in |-K_{V}|$ not passing
through $\pi^{-1}(p)$. In particular, the cubic curve $\pi(C)$
passes through all points of the set $\Sigma$ and does not pass
through the point $p$. Let $Y$ be a cone in ${\mathbb P}^4$ over
$\pi(C)$ with a vertex in a sufficiently general line in ${\mathbb
P}^4$. Then cubic hypersurface $Y$ passes through all points of
the set $\Sigma$ and does not pass through the point $p$.
\end{proof}

The following result is due to \cite{ChPa04}.

\begin{lemma}
\label{lemma:of-excluding-three-points} Let $\Delta\subset{\mathbb
P}^n$ be a subset and $p\in{\mathbb P}^n\setminus\Delta$ be a
point such that $\{p\cup\Delta\}\subset{\mathbb P}^n$ is not
contained in a linear subspace of dimension $r$. Then there is a
linear subspace $H\subset{\mathbb P}^n$ of dimension $r$ that
contains at least $r+1$ points of the set $\Delta$ but not the
point $p$.
\end{lemma}

\begin{proof}
We will prove the claim by induction on $n$. For $n=2$ the claim
is trivial. Suppose that $n>2$ and $r<n$. By assumption there are
$r+1$ points $\{q_1,\cdots, q_{r+1}\}\subset\Delta$ such that the
linear span $T$ of the points $q_{i}$ has dimension $r$. We may
assume $p\in T$, because otherwise we are done. Thus, there is a
point $q\in\Delta\setminus T$, because by assumption the subset
$\{p\cup\Delta\}\subset {\mathbb P}^n$ is not contained in a
linear subspace of dimension $r$. By induction there is a linear
subspace $S\subset T$ of dimension $r-1$ that contains $r$ points
among $\{q_1,\cdots, q_{r+1}\}$ but not $p$. Now consider a cone
$H\subset{\mathbb P}^n$ over $T$ with the vertex $q$. The cone $H$
is a linear subspace of dimension $r$ that contains at least $r+1$
points of the set $\Delta$ but not the point $p$.
\end{proof}

\begin{proposition}
\label{proposition:nodes-hyperplane} Let $\Gamma\subset{\mathbb
P}^4$ be a hyperplane such that ${\textsf{Sing}}(X)$ is contained
in $\Gamma$. Then the nodes of the quartic $X$ impose independent
linear conditions on cubic surfaces in $\Gamma\cong{\mathbb P}^3$
and on cubic hypersurfaces in ${\mathbb P}^4$.
\end{proposition}

\begin{proof}
Let $\Sigma\subsetneq{\textsf{Sing}}(X)$ be any subset and
$p\in{\textsf{Sing}}(X)\setminus\Sigma$ be a point. We must show
that there is a cubic surface in $\Gamma$ and a cubic hypersurface
in ${\mathbb P}^4$ passing through $\Sigma$ and not passing
through the point $p$. As in the proof of
Proposition~\ref{proposition:nodes-plane} it is enough to find a
cubic surface in $\Gamma$ that passes through all the points of
$\Sigma$ and does not pass through the point $p$. A sufficiently
general cone over such cubic surface gives a cubic hypersurface in
${\mathbb P}^4$ passing through all the points in the set $\Sigma$
and not passing through the point $p$. Without loss of generality
we may assume that $|{\textsf{Sing}}(X)|=|\Sigma|+1=9$.

Let $r\ge 2$ be the maximal number of points of the set $\Sigma$
that belong to a two-dimensional linear subspace $\Pi$ in $\Gamma$
together with $p$. Then $r\leq 7$ by
Proposition~\ref{proposition:nodes-in-general-position}. Let
$\Sigma=\{p_{1},\cdots,p_{8}\}$ and the first $r$ points of
$\Sigma$, i.e. the points $p_1,\cdots,p_r$, are contained in the
plane $\Pi$ together with $p$. Then the points $p$ and
$p_1,\cdots,p_r$ do not lie on a line, because otherwise we can
find a hyperplane in $\Gamma$ containing more than $r$ points of
the set $\Sigma$. We will prove the statement case by case.

Suppose $r=2$. Divide the set $\Sigma$ into three possibly
overlaping subsets such that each subset contains three points of
the set $\Sigma$ and their union is the whole set $\Sigma$. The
hyperplane in $\Gamma$ generated by each subset does not contain
$p$, because $r=2$. Hence, the product of these three hyperplanes
is the required cubic surface.

Suppose $r=3$. By Lemma~\ref{lemma:of-excluding-three-points}, we
can find three points of $\Sigma$ outside of $\Pi$, say
$p_{4},p_{5},p_{6}$, such that they generate the hyperplane in
$\Gamma$ not passing though $p$. Moreover, the four points
$\{p,p_1,p_2,p_3\}$ do not lie on one line by
Proposition~\ref{proposition:nodes-in-general-position}.
Therefore, there is a line passing through two points of the set
$\{p_1,p_2,p_3\}$, say through $p_{1}$ and $p_{2}$, and not
passing through the point $p$. Therefore, the product of the
hyperplane passing through the points $p_{4},p_{5},p_{6}$ and a
hyperplane passing through the points $p_{7},p_{1},p_{2}$ and a
sufficiently general hyperplane passing through the points $p_{3}$
and $p_{8}$ gives a cubic surface in $\Gamma\cong{\mathbb P}^3$
passing through all the points of the set $\Sigma$ and not passing
through the point $p$.

Suppose $r=4$. There are two lines in $\Pi$, say $L_{1}$ and
$L_{2}$, such that $L_{1}$ contains $p_{1}$ and $p_{2}$, the line
$L_{2}$ contains $p_{3}$ and $p_{4}$, and both lines do not pass
through $p$. Moreover, there are at most two points among
$\{p_{5},p_{6},p_{7},p_{8}\}$ that lie on a line passing through
the point $p$. Therefore, there are two points, say $p_{5}$ and
$p_{6}$, such that the line passing through the points $p_{5}$ and
$p_{6}$ does not pass through the point $p$. The product of two
hyperplanes passing through the lines $L_{1}$ and $L_{2}$ and two
points $p_{7}$ and $p_{8}$ respectivly and a sufficiently general
hyperplane passing through the points $p_{5}$ and $p_{6}$ gives
the required cubic surface.

Suppose $r=5$. There are two lines in $\Pi$, say $L_{1}$ and
$L_{2}$, such that $p\not\in L_{1}\cup L_{2}$ and $L_{1}\cup
L_{2}$ contains four points of $\Sigma\cap\Pi$, say the points
$p_{1},p_{2},p_{3}$ and $p_{4}$. The product of two hyperplane
passing through the lines $L_{1}$ and $L_{2}$ and two points
$p_{7}$ and $p_{8}$ respectivly and a sufficiently general
hyperplane passing through the points $p_{5}$ and $p_{6}$ gives a
cubic surface in $\Gamma$ passing through all the points of the
set $\Sigma$ and not passing through the point $p$.

Suppose $r=6$. Now we have six points of the set $\Sigma\cap\Pi$
and two points, say $p_{7}$ and $p_{8}$, of the set $\Sigma$
outside of $\Pi$. We can find a cubic curve $C$ on $\Pi$ that
passes through $\Sigma\cap\Pi$ and does not pass through $p$ by
Proposition~\ref{proposition:nodes-plane}. A sufficiently general
hyperplane in $\Gamma$ passing through the points $p_{7}$ and
$p_{8}$ meets the curve $C$ at three points. Let $q$ and
$q^{\prime}$ be two points among them and $O$ be an intersection
of the lines $<p_{7},q>$ and $<p_{8},q^{\prime}>$. Then the cubic
cone in $\Gamma$ over the curve $C$ with the vertex $O$ is a cubic
surface that passes through all the point of $\Sigma$ but not
through $p$.

Suppose $r=7$. We can find a cubic curve $C$ on $\Pi$ that passes
through the seven points of the set $\Sigma\cap\Pi$  and does not
pass through the point $p$ by
Proposition~\ref{proposition:nodes-plane}. The cone in
$\Gamma\cong{\mathbb P}^3$ over the cubic curve $C$ with the
vertex $p_{8}$ is a cubic surface that passes through $\Sigma$ but
not through $p$.
\end{proof}

\begin{proposition}
\label{proposition:nodes-general} The nodes of $X$ impose
independent linear conditions on cubics in ${\mathbb P}^4$.
\end{proposition}

\begin{proof}
We must show that for any subset
$\Sigma\subsetneq{\textsf{Sing}}(X)$ and a point
$p\in{\textsf{Sing}}(X)\setminus\Sigma$ there is a cubic
hypersurface in ${\mathbb P}^4$ passing through all the points of
$\Sigma$ and not passing through $p$. Without loss of generality
we may assume $|{\textsf{Sing}}(X)|=|\Sigma|+1=9$.

Let $r\ge 3$ be the maximal number of points in $\Sigma$ that
belongs to a hyperplane $\Xi\subset{\mathbb P}^4$ together with
$p$. We may assume $r\leq 7$ by
Proposition~\ref{proposition:nodes-hyperplane}. Let
$\Sigma=\{p_{1},\cdots,p_{8}\}$ and the first $r$ points of
$\Sigma$, i.e. the points $p_1,\cdots,p_r$, are contained in $\Xi$
together with $p$. Then the points $p$ and $p_1,\cdots,p_r$ do not
belong to a two-dimensional linear subspace in ${\mathbb P}^4$,
because otherwise we can find a hyperplane passing through $r+1$
points of $\Sigma$. We will prove the claim case by case.

Suppose $r=3$. Divide the set $\Sigma$ into three possibly
overlaping subsets such that each subset contains exactly four
points of the set $\Sigma$. The hyperplane generated by each
subset does not contain the point $p$, because $r=3$. The product
of these three hyperplanes is the required cubic hypersurface.

Suppose $r=4$. There are two lines $L_{1}$ and $L_{2}$ in $\Xi$
such that $L_{1}$ passes through  $p_{1}$ and $p_{2}$, the line
$L_{2}$ passes through $p_{3}$ and $p_{4}$, and both lines do not
pass through $p$. Moreover, there are at most two points of the
set $\{p_{5},p_{6},p_{7},p_{8}\}$ that lie on a line containing
$p$. Hence, there are two points, say $p_{5}$ and $p_{6}$, such
that the line passing through $p_{5}$ and $p_{6}$ does not pass
through $p$. The product of two sufficiently general hyperplanes
passing through the lines $L_{1}$ and $L_{2}$ and two points
$p_{7}$ and $p_{8}$ respectivly and a sufficiently general
hyperplane passing through the points $p_{5}$ and $p_{6}$ gives
the required cubic hypersurface in ${\mathbb P}^4$.

Suppose $r=5$. As in the previous case there are two lines $L_{1}$
and $L_{2}$ in $\Xi$ such that $L_{1}$ passes through the points
$p_{1}$ and $p_{2}$, line $L_{2}$ passes through the points
$p_{3}$ and $p_{4}$, and both lines do not pass through $p$. The
product of two general hyperplanes passing through the lines
$L_{1}$ and $L_{2}$ and through the points $p_{7}$ and $p_{8}$
respectivly and a sufficiently general hyperplane passing through
the points $p_{5}$ and $p_{6}$ gives a cubic hypersurface in
${\mathbb P}^4$ that passes through all the points of $\Sigma$ and
does not pass through the point $p$.

Suppose $r=6$. There are six points in $\Sigma\cap\Xi$ and two
points, say $p_{7}$ and $p_{8}$, of $\Sigma$ outside of the
hyperplane $\Xi$. There is a cubic surface $S\subset\Xi$ that
passes through the six points of $\Sigma\cap\Xi$ and does not pass
through $p$ by Proposition~\ref{proposition:nodes-hyperplane}. A
general two-dimensional linear subspace passing through the points
$p_{7}$ and $p_{8}$ meets $S$ at three different points. Choose
two points $q$ and $q^{\prime}$ among these intersection points.
Let $O$ be an intersection of the lines $<p_{7},q>$ and
$<p_{8},q^{\prime}>$. Now the required cubic hypersurface is a
cone in ${\mathbb P}^4$ over the cubic surface $S$ with the vertex
$O$.

Suppose $r=7$. We can find a cubic surface $S\subset\Xi$ that
passes through the seven points of the set $\Sigma\cap\Pi$ and
does not pass through $p$ by
Proposition~\ref{proposition:nodes-hyperplane}. The cone in
${\mathbb P}^4$ over the surface $S$ with the vertex $p_{8}$
passes through all the point of $\Sigma$ but not through $p$.
\end{proof}

Therefore, both Theorems~\ref{theorem:main} and
\ref{theorem:second} are proven. Apriori the same method can be
applied to any nodal hypersurface in ${\mathbb P}^4$. The
following result (cf. \cite{CiGe03}) is implied by \cite{We87} and
\cite{Di90}.

\begin{theorem}
\label{theorem:factorial-hypersurfaces} A nodal hypersurface
$V\subset{\mathbb P}^4$ is ${\mathbb Q}$-factorial when
$|{\textsf{Sing}}(V)|\le 2\deg(V)-4$.
\end{theorem}

The bound for nodes in
Theorem~\ref{theorem:factorial-hypersurfaces} is not sharp except
for hyperquadrics.

\section{The proof of Theorem~\ref{theorem:quartic-del-Pezzo}}
\label{section:proof-of-theorem-del-Pezzo}

Let $X\subset{\mathbb P}^4$ be a sufficiently general\footnote{A
complement to a Zariski closed subset in moduli.} quartic 3-fold
containing a smooth del Pezzo surface $S\subset{\mathbb P}^4$ of
degree $4$. Then the quartic $X$ can be given by the equation
$$
a_{2}(x,y,z,t,w)h_{2}(x,y,z,t,w)+b_{2}(x,y,z,t,w)g_{2}(x,y,z,t,w)=0\subset{\textsf{Proj}}({\mathbb C}[x,y,z,t,w]),%
$$
where $a_{2}$, $b_{2}$, $h_{2}$ and $g_{2}$ are homogeneous
polynomials of degree $2$ such that $S$ is defined by the
equations $h_{2}=g_{2}=0$. The quartic $X$ is nodal and
non-${\mathbb Q}$-factorial. Moreover, it has $16$ nodes given by
the equations $h_{2}=g_{2}=a_{2}=b_{2}=0$.

\begin{lemma}
\label{lemma:class-group} The divisor class group $\textsf{Cl}(X)$
is ${\mathbb Z}\oplus {\mathbb Z}$.
\end{lemma}

\begin{proof}
Let $f:U\to {\mathbb P}^4$ be a blow up of the surface $S$, $E$ be
an exceptional divisor of the birational map $f$ and
$H=f^{*}({\mathcal O}_{{\mathbb P}^4}(1))$. Then the pencil
$|2H-E|$ is free, because the surface $S$ is a complete
intersection of two quadrics in ${\mathbb P}^4$. In particular,
the divisor $2H-E$ is nef and the divisor $4H-E$ is ample. On the
other hand, the proper transform ${\tilde X}\subset U$ of the
quartic $X$ is rationally equivalent to the divisor $4H-E$. The
restriction $f\vert_{\tilde X}:{\tilde X}\to X$ is a small
resolution and ${\tilde X}$ is smooth. Therefore
$$
H^{2}({\tilde X}, {\mathbb Z})\cong H^{2}(U, {\mathbb Z})\cong {\mathbb Z}\oplus {\mathbb Z}%
$$
by the Lefschetz theorem (see \cite{AndFr59}, \cite{Bott59},
\cite{Mi63}, \cite{KuKur98}), which implies the claim of the
lemma.

The second way to prove the claim is to prove that the nodes of
$X$ impose $15$ independent linear conditions on cubic
hypersurfaces in ${\mathbb P}^4$, which implies the claim due to
\cite{We87}, \cite{Di90}, \cite{Cy01}. It is enough to prove that
the nodes of $X$ impose $15$ independent linear conditions on the
global sections of the sheaf ${\mathcal O}_{{\mathbb
P}^4}(3)\vert_{S}$ due to the surjectivity of the map
$H^{0}({\mathcal O}_{{\mathbb P}^4}(3))\to H^{0}({\mathcal O}_{{\mathbb P}^4}(3)\vert_{S})$. %
The latter can be deduced from \cite{Bes83} using
Proposition~\ref{proposition:nodes-in-general-position} and the
fact that $S$ is a blow up of ${\mathbb P}^{2}$ in $5$ points.
\end{proof}

The pencil generated by the quadrics $a_{2}=0$ and $b_{2}=0$ cuts
on the quartic $X$ the surface $S$ and a pencil ${\mathcal M}$
whose general element is a smooth del Pezzo surface of degree $4$.
Let $\tau:V\to X$ be a small resolution (see \cite{Ka88}) such
that the pencil ${\mathcal H}=\tau^{-1}({\mathcal M})$ is free,
i.e. $V={\textsf {Proj}}(\oplus_{i\ge 0}{\mathcal
O}_{X}(-S)^{\otimes i})$ and $\tau$ is a natural projection to
$X$. Then $V$ is smooth and projective,
${\textsf{Pic}}(V)={\mathbb Z}\oplus{\mathbb Z}$, and the pencil
${\mathcal H}$ gives a morphism $\tau:V\to{\mathbb P}^1$ whose
general fiber is a del Pezzo surface of degree $4$.

\begin{corollary}
\label{corollary:del-Pezzo-conic-bundle} The 3-fold $V$ is
bi\-ra\-ti\-o\-nal to a conic bundle (see \cite{Ma66},
\cite{Ma67}, \cite{Is72}, \cite{Is80a}, \cite{Al87},
\cite{Is96b}).
\end{corollary}

The generality in the choice of $X$ implies that $\tau$ is
standard in the sense of \cite{Al87}, i.e. every fiber of the
fibration $\tau$ is normal and $\textsf{Pic}(V)={\mathbb
Z}\oplus{\mathbb Z}$. The following result was proven in
\cite{Al87}.

\begin{theorem}
\label{theorem:Valera} Let $\gamma:Y\to{\mathbb P}^{1}$ be a
standard fibration of del Pezzo surfaces of degree $4$. Then $Y$
is non-rational if the topological Euler characteristic of the
3-fold $Y$ is different from $0$, $-8$ and $-4$.
\end{theorem}

Therefore, in order to prove
Theorem~\ref{theorem:quartic-del-Pezzo} we must calculate the
topological Euler characteristic of the 3-fold $V$. The following
result was proven in \cite{Cy01}.

\begin{theorem}
\label{theorem:Cynk} Let $W$ be a projective smooth 4-fold and
$Y\subset W$ be a reduced and reducible divisor such that the only
singularities of $Y$ are nodes. Let ${\tilde Y}$ be a small
resolution of $Y$. Suppose that
$$
h^{2}(\Omega^{1}_{W})=h^{3}(\Omega^{1}_{W}\otimes {\mathcal O}_{W}(-Y))=0%
$$
and the line bundle ${\mathcal O}_{W}(Y)$ is ample. Then
$h^{1}(\Omega^{1}_{\tilde Y})=h^{1}(\Omega^{1}_{W})+\delta$ and
$h^{2}(\Omega^{1}_{\tilde Y})$ equals to
$$
h^{0}(K_{W}\otimes {\mathcal O}_{W}(2Y)) %
+h^{3}({\mathcal O}_{W})-h^{0}(K_{W}\otimes {\mathcal O}_{W}(Y)) %
-h^{3}(\Omega^{1}_{W})-h^{4}(\Omega^{1}_{W}\otimes {\mathcal O}_{W}(-Y)) %
-|\textsf{Sing}(Y)|+\delta,
$$
where $\delta$ is the number of dependent equations that vanishing
at the nodes of $Y$ imposes on the global sections of the line
bundle $K_{W}\otimes {\mathcal O}_{W}(2Y)$, i.e. the defect of the
3-fold $Y$.
\end{theorem}

The topological Euler characteristic $\chi(V)$ of $V$ is
$6-2h^{2}(\Omega^{1}_{V})$. The twisted Euler exact sequence and
the Serre duality imply $h^{3}(\Omega^{1}_{{\mathbb
P}^4}\otimes{\mathcal O}_{{\mathbb P}^4}(-4))=0$ and
$h^{4}(\Omega^{1}_{{\mathbb P}^4}\otimes{\mathcal O}_{{\mathbb
P}^4}(-4))=5$. By Theorem~\ref{theorem:Cynk}
$$
h^{2}(\Omega^{1}_{V})=h^{0}({\mathcal O}_{{\mathbb P}^4}(3))-h^{3}(\Omega^{1}_{{\mathbb P}^4})-h^{4}(\Omega^{1}_{{\mathbb P}^4}\otimes{\mathcal O}_{{\mathbb P}^4}(-4))-|{\textsf{Sing}}(X)|+1=15%
$$
and $\chi(V)=-24$. The quartic $X$ is non-rational by
Theorem~\ref{theorem:Valera}. Hence,
Theorem~\ref{theorem:quartic-del-Pezzo} is proven.

\section{The proof of Theorem~\ref{theorem:quartic-plane}}
\label{section:proof-of-theorem-plane}

Let $X\subset{\mathbb P}^4$ be a very general\footnote{A
complement to a countable union of Zariski closed subsets in
moduli.} quartic 3-fold containing a plane $\Pi\subset{\mathbb
P}^4$. Then $X$ can be given by
$$
xh_{3}(x,y,z,t,w)+yg_{3}(x,y,z,t,w)=0\subset{\textsf{Proj}}({\mathbb C}[x,y,z,t,w]),%
$$
where $h_{3}$ and $g_{3}$ are homogeneous polynomials of degree
$3$, and the plane $\Pi$ is defined by the equations $x=y=0$. The
quartic $X$ is nodal, it has $9$ nodes given by
$x=y=h_{3}=g_{3}=0$.

\begin{lemma}
\label{lemma:class-group-generic} The divisor class group
$\textsf{Cl}(X)$ is ${\mathbb Z}\oplus {\mathbb Z}$.
\end{lemma}

\begin{proof}
The claim of the lemma is equivalent (see \cite{We87},
\cite{Di90}, \cite{Cy01}) to the following: the nodes of the
quartic $X$ impose $8$ independent linear conditions on cubic
hypersurfaces in ${\mathbb P}^4$, i.e. the defect of the quartic
$X$ is one (see Theorem~\ref{theorem:Cynk}). However, during the
proof of Theorem~\ref{theorem:main} we implicitly proved that any
$8$ nodes of the quartic $X$ impose $8$ independent linear
conditions on cubic hypersurfaces in ${\mathbb P}^4$. On the other
hand, the nodes of $X$ can not impose $9$ independent linear
conditions on cubic hypersurfaces in ${\mathbb P}^4$, because $X$
is obviously not ${\mathbb Q}$-factorial.
\end{proof}

To prove Theorem~\ref{theorem:quartic-plane} we will use the
degeneration technique (see \cite{Be77}, \cite{Tu80}, \cite{Cl82})
together with the following result in \cite{Ko95} and \cite{Ko96}.

\begin{theorem}
\label{theorem:Janos} Let $\xi:Y\to Z$ be a flat proper morphism
with irreducible and reduced geometric fibers. Then there are
countably many closed subvarieties $Z_{i}\subset Z$ such that for
an arbitrary closed point $s\in Z$ the fiber $\xi^{-1}(s)$ is
ruled if and only if $s\in\cup Z_{i}$.
\end{theorem}

Consider a sufficiently general quartic 3-fold $V\subset{\mathbb
P}^4$ given by the equation
$$
x{\bar h}_{3}(x,y,z,t,w)+y{\bar g}_{3}(x,y,z,t,w)=0\subset{\textsf{Proj}}({\mathbb C}[x,y,z,t,w])%
$$
such that
$$
{\bar h}_{3}(x,y,z,t,w)=xa_{2}(x,y,z,t,w)+yb_{2}(x,y,z,t,w)+f_{1}(z,t,w)h_{2}(z,t,w)%
$$
and
$$
{\bar g}_{3}(x,y,z,t,w)=xc_{2}(x,y,z,t,w)+yd_{2}(x,y,z,t,w)+f_{1}(z,t,w)g_{2}(z,t,w),%
$$
where $a_{2}$, $b_{2}$, $c_{2}$, $d_{2}$, $h_{2}$ and $g_{2}$ are
homogeneous polynomials of degree $2$, and $f_{1}$ is a
homogeneous polynomial of degree $1$. The quartic $X$ contains the
plane $\Pi$. The singularities of $X$ consist of $4$ nodes given
by the equations $x=y=h_{2}=g_{2}=0$ and a single double line
$L\subset\Pi$ given by the equations $x=y=f_{1}=0$.

\begin{remark}
\label{remark:degeneration} The resolution of singularities of $V$
has no global holomorphic forms and the Kodaira dimension of the
3-fold $V$ is $-\infty$, i.e.  the 3-fold $V$ is rationally
connected (see \cite{Ko96}). Hence, the rationality of $V$ is
equivalent to its ruledness. However, $V$ is a flat degeneration
of $X$. Thus, the non-rationality of $V$ implies the
non-rationality of $X$ by Theorem~\ref{theorem:Janos}.
\end{remark}

Therefore, to prove Theorem~\ref{theorem:quartic-plane} it is
enough to prove the non-rationality of the quartic $V\subset
{\mathbb P}^{4}$.

\begin{remark}
\label{remark:quartic-with-double-line} The non-rationality of a
sufficiently general quartic 3-fold with a double line was proven
in \cite{CoMu77a} and \cite{CoMu77b} using the method of
intermediate Jacobian (see \cite{ClGr72}, \cite{Be77},
\cite{Tu80}, \cite{Alb82}, \cite{Alz90}).
\end{remark}

Let $\pi:U\to {\mathbb P}^{4}$ be a blow up of the line $L\subset
{\mathbb P}^{4}$, $E$ be a $\pi$-exceptional divisor, ${\bar
V}\subset U$ be a proper transform of $V$. Then
$|\pi^{*}({\mathcal O}_{{\mathbb P}^4}(1))-E|$
is free and gives a ${\mathbb P}^{2}$-bundle $\psi:U\to{\mathbb P}^{2}$.%

\begin{lemma}
\label{lemma:conic-bundle} The 3-fold ${\bar V}$ is smooth in the
neighbourhood of the exceptional divisor $E$, the singularities of
${\bar V}$ consist of $4$ nodes, which are the images of the nodes
of $V$. For a point $x\in L$ the intersection $\pi^{-1}(x)\cap
{\bar V}\subset \pi^{-1}(x)\cong {\mathbb P}^{2}$ is a smooth
conic if $x$ is not a zero of $h_{2}$ or $g_{2}$ and a union of
two different lines otherwise, i.e. there are $4$ reducible fibers
of the morphism $\pi\vert_{E\cap {\bar V}}$.
\end{lemma}

\begin{proof}
Simple calculations.
\end{proof}

Let ${\bar \Pi}\subset U$ be a proper transform of $\Pi$. Then
$\psi(\bar \Pi)=O$ is a point. The restriction $\psi\vert_{\bar
V}:{\bar V}\to{\mathbb P}^{2}$ is a morphism whose fibers over
${\mathbb P}^{2}\setminus O$ are conics, and the fiber of
$\psi\vert_{\bar V}$ over the point $O$ is ${\bar \Pi}\subset
{\bar V}$.

\begin{lemma}
\label{lemma:standard-conic-bundle} Let $\gamma:W\to U$ be a blow
up of $\bar \Pi$, $G$ be a $\gamma$-exceptional divisor,
$\alpha:{\mathbb F}^{1}\to {\mathbb P}^{2}$ be a blow up of the
point $O$, and ${\tilde V}\subset W$ be a proper transform of
${\bar V}$. Then ${\tilde V}$ is a small resolution of the 3-fold
$\bar V$, the linear system $|\gamma^{*}(\pi^{*}({\mathcal
O}_{{\mathbb P}^4}(1))-E)-G|$ is free and gives a morphism
$\phi:W\to {\mathbb F}^{1}$ such that
$\psi\circ\gamma=\alpha\circ\phi$.
\end{lemma}

\begin{proof}
Simple calculations.
\end{proof}

\begin{lemma}
\label{lemma:Picard-group} The Picard group of the 3-fold $\tilde
V$ is ${\mathbb Z}\oplus {\mathbb Z}\oplus {\mathbb Z}$.
\end{lemma}

\begin{proof}
The divisor ${\tilde V}\subset W$ is rationally equivalent to a
divisor
$$
\gamma^{*}(\pi^{*}({\mathcal O}_{{\mathbb P}^4}(4))-2E)-G\sim
\gamma^{*}(\pi^{*}({\mathcal O}_{{\mathbb P}^4}(1))-E)-G+
\gamma^{*}(\pi^{*}({\mathcal O}_{{\mathbb P}^4}(1))-E)+
(\pi\circ\gamma)^{*}({\mathcal O}_{{\mathbb P}^4}(2))
$$
which is ample. Hence, $H^{2}({\tilde V}, {\mathbb Z})\cong
H^{2}(W, {\mathbb Z})$ by the Lefschetz theorem (see
\cite{AndFr59}, \cite{Bott59}, \cite{Mi63}, \cite{KuKur98}), which
implies the claim of the lemma.

Another way to prove the lemma is to prove that $\textsf{Cl}({\bar
V})$ is ${\mathbb Z}\oplus {\mathbb Z}\oplus {\mathbb Z}$. By
Theorem~\ref{theorem:Cynk} the latter is equivalent to the
following: the nodes of ${\bar V}$ impose $3$ independent linear
conditions on the global sections of the line bundle
$\pi^{*}({\mathcal O}_{{\mathbb P}^4}(3))-2E$. The latter is
implied by the following: the nodes of the 3-fold $V$ impose $3$
independent linear conditions on the hyperplanes in ${\mathbb
P}^4$ which is obvious.
\end{proof}

\begin{corollary}
\label{corollary-standard-conic-bundel} The restriction ${\tilde
\phi}=\phi\vert_{\tilde V}:{\tilde V}\to {\mathbb F}_{1}$ is a
standard conic bundle (see \cite{Sar80}).
\end{corollary}

Let $\Delta\subset {\mathbb F}_{1}$ be a degeneration divisor of
the standard conic bundle ${\tilde \phi}$. Then $\Delta$ is a
reduced divisor with at most simple normal crossing (see
\cite{Be77}, \cite{Tu80}, \cite{Sar80}, \cite{Sar82},
\cite{Sho84}, \cite{Co00}).

\begin{lemma}
\label{lemma-discriminant} Let $s_{\infty}$ be an exceptional
section of the ruled surface ${\mathbb F}_{1}$ and $l$ be a fiber
of the natural projection of the surface ${\mathbb F}_{1}$ to
${\mathbb P}^{1}$. Then $\Delta\sim 5s_{\infty}+8l$ and
$2K_{{\mathbb F}_{1}}+\Delta\sim s_{\infty}+2l$.
\end{lemma}

\begin{proof}
Let $\Delta\sim as_{\infty}+bl$ for some integer $a$ and $b$.
Consider a sufficiently general divisor $H$ in the linear system
$|{\tilde \phi}^{*}(l)|$ and the surface ${\tilde
\Pi}=\psi^{-1}(s_{\infty})$. Then $H$ and ${\tilde \Pi}$ are
smooth. Indeed, $H$ is smooth by the Bertini theorem, and the
surface $\tilde \Pi$ is smooth, because $\gamma\vert_{\tilde
\Pi}:{\tilde \Pi}\to {\bar \Pi}\cong\Pi$ is a blow up of the four
points on $\Pi\cong {\mathbb P}^{2}$ given by the equations
$h_{2}=g_{2}=0$. The birational map $\gamma\vert_{\tilde \Pi}$
resolves the base points of the pencil generated by the conics
$h_{2}=0$ and $g_{2}=0$, which induces the restriction morphism
$\phi\vert_{\tilde \Pi}$. The surface $H$ is a cubic surface,
whose image on the quartic $V$ is a cubic surface residual to the
plane $\Pi$. Hence, $K_{H}^{2}=3$ and $K_{{\tilde \Pi}}^{2}=5$.
Thus, $\Delta\cdot l=5$ and $\Delta\cdot s_{\infty}=3$.
\end{proof}

The following result in \cite{Sho84} is a special case of a
conjectural rationality criterion of a standard three-dimensional
conic bundle (see \cite{Is87}, \cite{Is91}, \cite{Is92},
\cite{Is96a}).

\begin{theorem}
\label{theorem:Slava} Let $\xi:Y\to Z$ be a conic bundle, $D$ be a
degeneration divisor of $\xi$, where $Y$ is a smooth 3-fold,
$\textsf{Pic}(Y/Z)={\mathbb Z}$ and $Z$ is either a ruled surface
${\mathbb F}_{r}$ or ${\mathbb P}^{2}$. Then the rationality of
the 3-fold $Y$ implies $|2K_{Z}+D|=\emptyset$.
\end{theorem}

Therefore, the 3-fold ${\tilde V}$ is non-rational by
Theorem~\ref{theorem:Slava}. Hence,
Theorem~\ref{theorem:quartic-plane} is proven.

\end{document}